\documentclass[12pt]{article}

\usepackage{amsmath}
\usepackage{amssymb}
\usepackage[dvips]{graphicx}

\begin{document}
\title{{\bf
{\normalsize The Ma-Qiu index and the Nakanishi index for a fibered knot are equal, 
and $\omega$-solvability}}
\footnotetext[0]{%
2020 {\it Mathematics Subject Classification}: 20F12, 20F16, 57K10, 57K31.\\
{\it Keywords}:
knot group; commutator subgroup; MQ index; Nakanishi index; 
$\omega$-solvale group.}}
\author{{\footnotesize Teruhisa KADOKAMI}}
\date{{\footnotesize February 19, 2023}}
\maketitle

\newcommand{\bsquare}{\hbox{\rule{6pt}{6pt}}}
\newcommand{\qed}{\hbox{\rule[-2pt]{3pt}{6pt}}}

\newtheorem{df}{Definition}[section]
\newtheorem{lm}[df]{Lemma}
\newtheorem{theo}[df]{Theorem}
\newtheorem{re}[df]{Remark}
\newtheorem{pr}[df]{Proposition}
\newtheorem{ex}[df]{Example}
\newtheorem{co}[df]{Corollary}
\newtheorem{cl}[df]{Claim}
\newtheorem{qu}[df]{Question}
\newtheorem{pb}[df]{Problem}

\makeatletter
\renewcommand{\theequation}{%
\thesection.\arabic{equation}}
\@addtoreset{equation}{section}
\makeatother

\begin{abstract}
{\footnotesize
\setlength{\baselineskip}{10pt}
\setlength{\oddsidemargin}{0.25in}
\setlength{\evensidemargin}{0.25in}
For a knot $K$ in $S^3$, let $G(K)$ be the knot group of $K$, 
$a(K)$ the Ma-Qiu index (the MQ index, for short), 
which is the minimal number of normal generators of 
the commutator subgroup of $G(K)$, and
$m(K)$ the Nakanishi index of $K$, 
which is the minimal number of generators of 
the Alexander module of $K$.
We generalize the notions for a pair of a group $G$ and its normal sugroup $N$, 
and we denote them by $a(G, N)$ and $m(G, N)$ respectively.
Then it is easy to see $m(G, N)\le a(G, N)$ in general.
We also introduce a notion ``$\omega$-solvability" for a group 
that the intersection of all higher commutator subgroups is trivial.
Our main theorem is that if $N$ is $\omega$-solvable, then we have $m(G, N)=a(G, N)$.
As corollaries, for a fibered knot $K$, we have $m(K)=a(K)$, and
we could determine the MQ indices of prime knots up to $9$ crossings completely.
}

\end{abstract}

\section{Introduction}\label{sec:intro}
For a group $G$ and its normal subgroup $N$, 
$y_1, \ldots, y_d\in N$ (possibly $d$ is infinite)
are {\it normal generators} of $N$ if any element $x$ in $N$ is 
a product of conjugates of them:
$$x=\prod_{j=1}^r
g_j^{-1}y_{i_j}^{\varepsilon_j}g_j\quad
(g_j\in G, \varepsilon_j=\pm 1, i_j\in \{1, \ldots, d\}; j=1, \ldots, r)$$
where $r$ is finite.
Then we denote by 
$N=\langle \! \langle y_1, \ldots, y_d\rangle \! \rangle$, and
the minimal number of such $d$ is called
the {\it Ma-Qiu index} (the {\it MQ index}, for short) of the group pair $(G, N)$,
which is denoted by $a(G, N)$.
It is easy to see that 
$a(G, N)=0$ if and only if $N$ is trivial (i.e.\ $N=\{1_G\}$).
J.~Ma and R.~Qiu \cite{MQ} originally defined the MQ index
for the case that $G$ is the knot group of a knot $K$ in $S^3$
and $N$ is the first commutator subgroup of $G$ (i.e.\ $N=[G, G]$),
and it is denoted by $a(K)$.
Then $a(K)=0$ if and only if $K$ is the trivial knot by the Dehn lemma
as pointed out in \cite{MQ}.
Their first main theorem is that the MQ index is a lower bound of 
the unknotting number (cf.\ Proposition \ref{pr:known} (1)).
For the proof, they used a {\it modified Wirtinger presentation}.
By the result, for the Kinoshita-Terasaka knot $K_{KT}$ \cite{KT}
and the Conway knot $K_C$ \cite{Co}, since they both have 
the unknotting number one \cite{BM} and they are both non-trivial, 
we have $a(K_{KT})=a(K_C)=1$.
Ma and Qiu also asked in \cite{MQ} about additivity of the MQ index
under connected sum (cf.\ Proposition \ref{pr:known} (2)).
The author and Yang \cite{KY} provided
counterexamples for additivity (``$1+1=1$")
of the MQ index (cf.\ Proposition \ref{pr:known} (5)).
We also showed that the MQ index is a lower bound
of both the rank$-1$ and the tunnel number
(cf.\ Proposition \ref{pr:known} (3) and (4)).
Note that Ma and Qiu did not use the term ``MQ index", 
but the author and Yang \cite{KY} firstly used the term.

\medskip

We regard $N/[N, N]$ as a group homology of $N$
where $[N, N]$ is the first commutator subgroup of $N$.
A multiplication in $N/[N, N]$ is an addition in the homology.
Moreover in the present case, 
a conjugation of $n\in N$ by $g\in G$ induces
a multiplication to $[n]\in N/[N, N]$ by $\overline{g} \in G/N$
where $[n]$ and $\overline{g}$ are the induced elements from $n$ and $g$ respectively
(i.e.\ $\left[g^{-1}ng\right]=[n]\cdot \overline{g}$\ 
($g\in G;\ n\in N$)).
Thus $N/[N, N]$ has a natural right $\mathbb{Z}[G/N]$-module structure.
We denote it by 
$$H(G, N)=H_1(N; \mathbb{Z}[G/N]),$$
and call it the {\it group homology} of $(G, N)$.
The {\it Nakanishi index} of the group pair $(G, N)$,
denoted by $m(G, N)$,
as the minimal number of generators of 
$H(G, N)$ as a right $\mathbb{Z}[G/N]$-module.
It is easy to see that $m(G, N)=0$ if and only if $H(G, N)=0$, and
$m(G, N)=1$ if and only if $H(G, N)$ is 
a non-trivial cyclic $\mathbb{Z}[G/N]$-module.
The original Nakanishi index, defined by Y.~Nakanishi \cite{Na},
is for the case that
$G$ is the knot group of a knot $K$ in $S^3$
and $N$ is the first commutator subgroup of $G$,
the index is the minimal size of the Alexander matrix of $K$,
and it is denoted by $m(K)$.
A.~Kawauchi \cite{Kw1} gave an alternative definition 
that $m(K)$ is minimal number of generators of the Alexander module of $K$.
Then $m(K)=0$ if and only if 
the Alexander polynomial of $K$, denoted by ${\mit \Delta}_K(t)$, is trivial
(i.e.\ ${\mit \Delta}_K(t)=1$)
(cf.\ \cite[Theorem 4.9.1]{Ne}).
It is easy to see that
\begin{equation}\label{eq:m&a}
m(G, N)\le a(G, N)
\end{equation}
because a set of normal generators realizing $a(G, N)$
also generates $H(G, N)$ via the projections
$G\to G/N$ and $N\to N/[N, N]$.
The equality of (\ref{eq:m&a}) does not hold in general.
For the Kinoshita-Terasaka knot $K_{KT}$ \cite{KT}
and the Conway knot $K_C$ \cite{Co}, 
we have $m(K_{KT})=m(K_C)=0$ and $a(K_{KT})=a(K_C)=1$.

\medskip

For a group $H$, we denote the $n$-th commutator subgroup of $H$ by $D^{(n)}(H)$.
In particular, $D^{(0)}(H)=H$ and $D^{(1)}(H)=[H, H]$.
Then $H$ is solvable if there exists a finite $n$ such that $D^{(n)}(H)=\{1_H\}$.
We define that $H$ is {\it $\omega$-solvable} 
if the intersection of all $D^{(n)}(H)$ is $\{1_H\}$.
It is easy to see that a solvable group is $\omega$-solvable.
A free group with rank greater than one is not solvable, but $\omega$-solvable
(see Section \ref{sec:known}).

\medskip

The following is our main theorem:

\begin{theo}\label{th:main}
For a group pair $(G, N)$, 
if $N$ is $\omega$-solvable, then we have $m(G, N)=a(G, N)$.
\end{theo}

For a fibered knot $K$, since the commutator subgroup of $G(K)$ is 
a free group with rank greater than one, 
which is the fundamental group of the fiber surface,
$G(K)$ is $\omega$-solvable.
Therefore we have:

\begin{co}\label{co:fiber}
Let $K$ be a fibered knot in $S^3$.
Then we have $m(K)=a(K)$.
\end{co}

This is the affirmative answer for Question 5.7 in \cite{KY}.
As corollaries, 
(1) one of the main theorems in \cite{KY}, Theorem 1.3, 
is obtained immediately 
(cf.\ Proposition \ref{pr:known} (5)), and 
(2) for prime knots up to 9 crossings, the MQ index is determined completely
(see Section \ref{sec:table}).

\section{Preliminaries}\label{sec:known}

For a knot $K$ in $S^3$,
let $u(K)$ be the {\it unknotting number} of $K$,
$t(K)$ the {\it tunnel number} of $K$, and
$r(K)$ the {\it rank} (of the knot group) of $K$.
Ma and Qiu \cite{MQ}, and
the author and Z.~Yang \cite{KY}
showed the following :
\begin{pr}\label{pr:known}
Let $K$, $K_1$ and $K_2$ be knots in $S^3$.
\begin{enumerate}
\item[(1)]
{\rm \cite{MQ}}
$m(K)\le a(K)\le u(K).$

\item[(2)]
{\rm \cite{MQ}}
$\max\{a(K_1), a(K_2)\}\le a(K_1\sharp K_2)\le a(K_1)+a(K_2)$.

\item[(3)]
{\rm \cite[Theorem 1.1]{KY}}
$a(K)\le r(K)-1$.

\item[(4)]
{\rm \cite[Corollary 1.2]{KY}}
$m(K)\le a(K)\le \min \{r(K)-1, u(K)\} \le \min \{t(K), u(K)\}.$

\item[(5)]
{\rm \cite[Theorem 1.3]{KY}}
For odd integers $p$ and $q$ with $|p|, |q|\ge 3$,
let $K_{p, q}$ be the connected sum of
the $(2, p)$-torus knot and the $(2, q)$-torus knot.
Then $m(K_{p, q})=a(K_{p, q})=1$ or $2$.
Moreover they are $1$ if and only if $\gcd(p, q)=1$.

\end{enumerate}
\end{pr}

Here we give the definition of ``$\omega$-solvability".
For a group $H$, the $n$-th commutator subgroup of $H$, 
denoted by $D^{(n)}(H)$, is defined inductively by (C1) and (C2) as follows:

\medskip

\noindent
(C1) : $D^{(0)}(H)=H$ (and $D^{(1)}(H)=[H, H]$),

\medskip

\noindent
(C2) : $D^{(n+1)}(H)=[D^{(n)}(H), D^{(n)}(H)]$,

\medskip

\noindent
where for subgroups $H_1$ and $H_2$ of $H$, 
$[H_1, H_2]$ is the commutator subgroup generated by $H_1$ and $H_2$.
We note that $D^{(n)}(H)$ is a characteristic subgroup 
(i.e.\ a special case of a normal subgroup, which is stable under 
any automorphism of $H$), and
\begin{equation}\label{eq:inclusion}
D^{(0)}(H)=H\supset D^{(1)}(H)=[H, H]\supset
\cdots \supset D^{(n)}(H)\supset D^{(n+1)}(H)\supset \cdots
\end{equation}
Then we denote by
\begin{equation}\label{eq:omega}
D^{(\omega)}(H)=\bigcap_{n=0}^{\infty}D^{(n)}(H)
\end{equation}
and we define that $H$ is {\it $\omega$-solvable} 
if $D^{(\omega)}(H)=\{1_H\}$.
It is easy to see that a $\omega$-solvable group is a characteristic subgroup, 
and a solvable group is $\omega$-solvable.
A free group with rank greater than one is not solvable, but $\omega$-solvable.

\section{Proof of Theorem \ref{th:main}}\label{sec:proof}

We take $z_1, \ldots, z_m\in H(G, N)$ such that 
$z_1, \ldots, z_m$ generate $H(G, N)$ as a right $\mathbb{Z}[G/H]$-module
realizing $m=m(G, N)$, and
$y_1, \ldots, y_m\in N$ such that $y_j$ is a lift of $z_j$\ $(j=1, \ldots, m)$
by the natural projection $p : N\to N/[N, N]$ 
(i.e.\ $p(y_j)=z_j$\ $(j=1, \ldots, m)$).
Then we set $Y=\{y_1, \ldots, y_m\}$ and
$\langle \! \langle y_1, \ldots, y_m\rangle \! \rangle
=\langle \! \langle Y\rangle \! \rangle
\subset N$.
By (\ref{eq:m&a}), 
we will only show $\langle \! \langle Y\rangle \! \rangle=N$.
Let $q_i : D^{(i)}(N)\to N/\langle \! \langle Y\rangle \! \rangle$\ 
$(i=0, 1, \ldots)$ 
be a homomorphism induced by the natural inclusion.
Let $\mathrm{Im}(q_i)$ be the image of $q_i$ in $N/\langle \! \langle Y\rangle \! \rangle$.
Then it is easy to see $\mathrm{Im}(q_i)\supset \mathrm{Im}(q_{i+1})$
by (\ref{eq:inclusion}).
We will show $\mathrm{Im}(q_i)\subset \mathrm{Im}(q_{i+1})$
(as a result, $\mathrm{Im}(q_i)=\mathrm{Im}(q_{i+1})$) 
for all $i$ ($i=0, 1, \ldots$).

\bigskip

We show the following key lemma:

\begin{lm}\label{lm:key}
For any $a_i, b_i\in D^{(i)}(N)$, 
there exist $a_{i+1}, b_{i+1}\in D^{(i+1)}(N)$ and
$c_i, d_i\in \langle \! \langle Y\rangle \! \rangle$
such that $a_i=c_ia_{i+1}, b_i=d_ib_{i+1}$ and
$$[a_i,b_i]\equiv [a_{i+1},b_{i+1}]\quad
(\mathrm{mod}\ \langle \! \langle Y\rangle \! \rangle)$$
for all $i$ $(i=0, 1, \ldots)$.
\end{lm}

\noindent
{\bf Proof}\ \ 
Suppose that there exist $c_i, d_i\in \langle \! \langle Y\rangle \! \rangle$
such that $a_i=c_ia_{i+1}, b_i=d_ib_{i+1}$.
Then we have:
$$
\begin{array}{rcl}
[a_i,b_i] & = & a_i^{-1}b_i^{-1}a_ib_i\medskip\\
& = & (c_ia_{i+1})^{-1}(d_ib_{i+1})^{-1}(c_ia_{i+1})(d_ib_{i+1})\medskip\\
& = & a_{i+1}^{-1} c_i^{-1} b_{i+1}^{-1} d_i^{-1} c_ia_{i+1} d_ib_{i+1}\medskip\\
& = & a_{i+1}^{-1} c_i^{-1} b_{i+1}^{-1} d_i^{-1} c_i(a_{i+1} d_ia_{i+1}^{-1})
a_{i+1}b_{i+1}\medskip\\
& = & a_{i+1}^{-1} c_i^{-1} (b_{i+1}^{-1} (d_i^{-1} c_i(a_{i+1} d_ia_{i+1}^{-1}))b_{i+1})
b_{i+1}^{-1}a_{i+1}b_{i+1}\medskip\\
& = & (a_{i+1}^{-1} (c_i^{-1} (b_{i+1}^{-1} (d_i^{-1} c_i(a_{i+1} d_ia_{i+1}^{-1}))b_{i+1}))a_{i+1})
[a_{i+1},b_{i+1}]\medskip\\
& \equiv & [a_{i+1},b_{i+1}]\quad
(\mathrm{mod}\ \langle \! \langle Y\rangle \! \rangle).
\end{array}
$$

We show the lemma by induction on $i$.

\medskip

\noindent
(1)\ The case $i=0$.

\medskip

By the assumption on $p : N\to N/[N, N]$ 
and $Y=\{y_1, \ldots, y_m\}$, 
for any $a_0, b_0\in D^{(0)}(N)=N$, 
there exist $a_1, b_1\in D^{(1)}(N)=[N, N]$ 
and $c_0, d_0\in \langle \! \langle Y\rangle \! \rangle$
such that $a_0=c_0a_1, b_0=d_0b_1$
(i.e.\ $p(a_0)=p(c_0), p(b_0)=p(d_0)$).
By the calculation above, we have
$[a_0,b_0]\equiv [a_1,b_1]\quad
(\mathrm{mod}\ \langle \! \langle Y\rangle \! \rangle)$.

\bigskip

\noindent
(2)\ Suppose the case $i$\ $(i\ge 0)$ is true.
We show the case $(i+1)$.

\medskip

By the assumption on $i$, 
for any $a_i, b_i\in D^{(i)}(N)$, 
there exist $a_{i+1}, b_{i+1}\in D^{(i+1)}(N)$ 
and $c_i, d_i\in \langle \! \langle Y\rangle \! \rangle$
such that $a_i=c_ia_{i+1}, b_i=d_ib_{i+1}$.
By the calculation above, we have
$[a_i,b_i]\equiv [a_{i+1},b_{i+1}]\quad
(\mathrm{mod}\ \langle \! \langle Y\rangle \! \rangle)$.
Since any element of $D^{(i+1)}(N)$ is a product 
of the form $[a_i,b_i]$\ $a_i, b_i\in D^{(i)}(N)$, we have the result.
\qed

\bigskip

Now we go back the proof of Theorem \ref{th:main}.
Lemma \ref{lm:key} shows $\mathrm{Im}(q_i)\subset \mathrm{Im}(q_{i+1})$.
Since $\mathrm{Im}(q_i)\subset \mathrm{Im}(q_{i+1})$ for all $i$ ($i=0, 1, \ldots$),
$q_0$ is surjective, and 
$D^{(\omega)}(N)=\{1_N\}$ (cf.\ (\ref{eq:omega})) by the assumption,
we have that $N/\langle \! \langle Y\rangle \! \rangle$ is trivial.
Hence we have $N=\langle \! \langle Y\rangle \! \rangle$.
This completes the proof.
\qed

\bigskip

As corollaries of Theorem \ref{th:main}, 
we have Corollary \ref{co:fiber}, 
which is the affirmative answer for Question 5.7 in \cite{KY}, 
and Proposition \ref{pr:known} (5) 
(\cite[Theorem 1.3]{KY}) immediately.

\section{The MQ indices for prime knots up to $10$ crossings}\label{sec:table}

By Proposition \ref{pr:known} (4) and \cite{Kw2, Kw3, KW, MSY}, 
the MQ index for a prime knot up to $10$ crossings is 1 or 2.
By Corollary \ref{co:fiber}, we could determine that
the MQ indices of the following 26 knots are all one:
$$8_{16}, 9_{29}, 9_{32}, 10_{62}, 10_{64}, 10_{79}, 10_{81}, 10_{85}, 10_{89}, 
10_{94}, 10_{96}, 10_{100}, 10_{105}, 10_{106}, $$
$$10_{109}, 10_{110}, 10_{112}, 
10_{116}, 10_{148}, 10_{149}, 10_{150}, 10_{151}, 10_{152}, 10_{153}, 10_{154}, 
10_{158}.$$
On the other hand, 
the MQ indices of the following 19 knots are still open:
$$10_{65}, 10_{66}, 10_{67}, 10_{68}, 10_{80}, 10_{83}, 
10_{86}, 10_{87}, 10_{90}, 10_{92}, 10_{93}, $$
$$10_{97}, 10_{108}, 10_{111}, 
10_{117}, 10_{120}, 10_{121}, 10_{163}, 10_{166}.$$

\begin{co}\label{co:table}
Let $K$ be a prime knot up to $10$ crossings other than $19$ knots above.
Then we have $m(K)=a(K)$.
\end{co}

\section{Final remarks}\label{sec:remark}

In \cite{JKRS}, the MQ index is discussed as a lower bound 
of the unknotting number of a 2-knot.
In the present paper, we clarified that 
a group theoretic property ``$\omega$-solvability" 
induces the equality of the MQ index and the Nakanishi index.
Now we ask the opposite direction:

\begin{qu}\label{qu:opposite}
For a group pair $(G, N)$, 
if $m(G, N)=a(G, N)$, then is $N$ always $\omega$-solvable?
\end{qu}

A contrast notion of $\omega$-solvability may be perfectivity.
Actually, the Alexander polynomial of a non-trivial knot $K$ is equal to $1$ 
if and only if $[G(K), G(K)]$ is perfect, and then we have $m(K)=0$.
On the other hand, a knot $K$ is non-trivial if and only if $a(K)\ge 1$.
A group $H$ is {\it finite order perfect}/{\it higher order perfect}/{\it order $k$ perfect}
if $D^{(k)}(H)=D^{(k+1)}(H)\ne \{1_H\}$ for finite $k$.
We raise the following question:

\begin{qu}\label{qu:perfect}
For a group pair $(G, N)$, 
$m(G, N)<a(G, N)$ if and only if $N$ is finite order perfect?
\end{qu}

As the final remark of the paper, 
since we do not ask that $N$ should include $[G, G]$, 
Theorem \ref{th:main} can be applied for the case 
that $G/N$ is non-commutative.

\bigskip

{\noindent {\bf Acknowledgements}}\
The author is supported by JSPS Grant-in-Aid for Scientific Research(C) 21K03245.
The author thanks to Tetsuya Ito who informed him that 
the notion ``$\omega$-solvability" is nothing but the notion ``residually solvability".

{\footnotesize
Teruhisa KADOKAMI\par 
School of Mechanical Engineering,\par
College of Science and Engineering,\par
Kanazawa University,\par
Kakuma-machi, Kanazawa, Ishikawa, 920-1192, Japan\par
{\tt kadokami@se.kanazawa-u.ac.jp}\par
}

\end{document}